\documentclass[a4paper,11pt]{article}

\usepackage{authblk}
\usepackage{parskip}
\usepackage{mathtools}
\usepackage{amssymb}
\usepackage{amsmath}
\usepackage{latexsym}
\usepackage{mathrsfs}  
\usepackage{bbm}       
\usepackage{amsthm}    
\usepackage{relsize}   
\usepackage{graphicx}
\usepackage{thmtools}  
\usepackage{mathrsfs}  
\usepackage{paralist}  
\usepackage{lipsum}    
\usepackage[all]{xy}   

\numberwithin{equation}{section}

\usepackage{titlesec}   

\titleformat{\section}[block]{\bfseries\filcenter}
{{\upshape\thesection\enspace}}{.5em}{}

\titleformat{\subsection}[block]{\filcenter}
{{\upshape\thesubsection\enspace}}{.5em}{} 

\usepackage{enumitem}  
\setlist{nosep}  
\setitemize[0]{leftmargin=*}
\setenumerate[0]{leftmargin=*}
\setenumerate[1]{label={(\arabic*)}} 

\newcommand{\N}{\mathbb{N}}     
\newcommand{\R}{\mathbb{R}}     
\newcommand{\C}{\mathbb{C}}     
\newcommand{\Prob}{\mathbb{P}}  
\newcommand{\Exp}{\mathbb{E}}   
\newcommand{\goth}[1]{\mathfrak{#1}} 
\newcommand{\inner}[2]{\left\langle #1 \, , \, #2 \right\rangle} 
\newcommand{\norm}[1]{\left|\left|#1\right|\right|}              
\newcommand{\triplet}[3]{\left( #1, #2, #3 \right) }             
\newcommand{\ProbSpace}{\triplet{\Omega}{\mathscr{F}}{\Prob}}    
\newcommand{\abs}[1]{\left| #1 \right|}                          
\renewcommand{\qedsymbol}{$\square$}                       

\newcommand{\defeq}{\mathrel{\mathop:}=}                         
\newcommand\restr[2]{{
  \left.\kern-\nulldelimiterspace 
  #1 
  \vphantom{\big|} 
  \right|_{#2} 
  }}
  
\theoremstyle{plain} 
\newtheorem{theo}{Theorem}[section]    

\newtheorem{coro}[theo]{Corollary}
\newtheorem{lemm}[theo]{Lemma}
\newtheorem{assu}[theo]{Assumption}
\newtheorem{rema}[theo]{Remark}

\theoremstyle{definition} 
\newtheorem{defi}[theo]{Definition}
\newtheorem{exam}[theo]{Example}

\declaretheoremstyle[%
  spaceabove=-5pt,%
  spacebelow=6pt,%
  headfont=\normalfont\itshape,%
  postheadspace=1em,%
  qed=\qedsymbol%
]{mystyle} 
\declaretheorem[name={Proof},style=mystyle,unnumbered,
]{prf}

 \title{\Large Regularization of Cylindrical Processes In Locally Convex Spaces}
\author{C. A. Fonseca-Mora }
\affil{  Escuela de Matem\'{a}tica, Universidad de Costa Rica, San Jos\'{e}, \\
Cod. 11501-2060, Costa Rica. \\
\noindent E-mail:  christianandres.fonseca@ucr.ac.cr }

\date{}

\begin{document}



 \maketitle

\abstract{Let $\Phi$ be a locally convex space and let $\Phi'$ denote its strong dual. In this paper we introduce sufficient conditions for the existence of a continuous or a c\`{a}dl\`{a}g $\Phi'$-valued version to a cylindrical process defined on $\Phi$. Our result generalizes many other known results on the literature and  their different connections will be discussed. As an application, we provide sufficient conditions for the existence of a $\Phi'$-valued c\`{a}dl\`{a}g L\'{e}vy process version to a given cylindrical L\'{e}vy process in $\Phi'$. }

\smallskip



\section{Introduction}

An important topic within the theory of probability on infinite dimensional linear spaces is the study of cylindrical probability measures, that is, of finitely additive set functions that have ``projections'' on finite dimensional spaces that are probability measures. One of the main problems of study of such objects is to determine when a given cylindrical probability measure extend to a Radon probability measure. In \cite{Sazonov:1958}, Sazonov shown that within the context of a Hilbert space, a necessary and sufficient condition for a Radon extension of a cylindrical measure is the continuity of its Fourier transform with respect to a topology that is weaker than the norm topology. In a similar way, in \cite{Minlos:1963} Minlos discovered that in the dual of a countably Hilbertian nuclear space, the same necessary and sufficient condition holds but this time it is only required continuity on the nuclear topology. It was pointed out by Kolmogorov \cite{Kolmogorov} that in duals of countably Hilbertian spaces continuity of Fourier transforms with respect to what is nowadays known as Sazonov's topology is sufficient for the Radon extension. If the underlaying space is nuclear, the Sazonov topology coincides with the nuclear topology and hence both Sazonov and Minlos discoveries are part of a more general result. A further extension to more general classes of locally convex spaces was carried out later by  Badrikian \cite{Badrikian}, Bourbaki \cite{Bourbaki}, and Schwartz \cite{SchwartzRM}. Again, continuity of Fourier transform with respect to the Sazonov topology was shown to be sufficient for a cylindrical measure to have a Radon extension. This result if often known as the \emph{Sazonov-Minlos theorem}.

Closely related to the concept of cylindrical probability measures is that of cylindrical random variables. Such objects are generalised random variables defined on an infinite dimensional linear space and whose laws define only a cylindrical probability measure. Conservely, to every cylindrical probability measure  there corresponds a ``canonical'' cylindrical random variable defined on some probability space. 

Parallel to the problem of  finding conditions for the extension of a cylindrical probability measure to a probability measure in the dual $\Phi'$ of a locally convex space $\Phi$, is the problem of finding a $\Phi'$-valued random variable $Y$ that is a version of a given  cylindrical random variable $X$ in $\Phi'$, where here ``version'' means that for every $\phi \in \Phi$, the evaluation $X(\phi)$ of $X$ in $\phi$ is a version of the real-valued random variable $\inner{Y}{\phi}$ corresponding to $Y$ through duality. In \cite{ItoNawata:1983}, 
It\^{o} and Nawata introduced their \emph{regularization theorem} that states that under the assumption that $\Phi$ is a locally convex space whose topology is generated by a family of separable Hilbertian seminorns (i.e. a multi-Hilbertian space),  the continuity of the Fourier transform of the cylindrical random variable on the Sazonov topology on $\Phi$ is a sufficient condition for the existence of a $\Phi'$-valued version. 

Now, if $(X_{t}: t \geq 0)$ is a cylindrical process in $\Phi'$, i.e. a collection of cylindrical random variables in $\Phi'$, is such that for each $t \geq 0$ the Fourier transform of $X_{t}$ is continuous in the Sazonov topology on $\Phi$, the regularization theorem of It\^{o} and Nawata provides the existence of a $\Phi'$-valued process $(Y_{t}: t \geq 0)$ such that the random variable $Y_{t}$ is a version of the cylindrical random variable $X_{t}$ for each $t \geq 0$. However, no statement can be made about the path properties (i.e. the regularity) of the process $(Y_{t}: t \geq 0)$. For example, even if for each $\phi \in \Phi$ the real-valued process $(X_{t}(\phi): t \geq 0)$ has continuous paths, using the regularization theorem alone it is not in general possible to conclude that  $(Y_{t}: t \geq 0)$ has continuous paths in $\Phi'$. 

For the reasons given above, in  \cite{Mitoma:1983} and under the assumption that $\Phi$ is a nuclear Fr\'{e}chet space, Mitoma shown that if $(X_{t}: t \geq 0)$ is a $\Phi'$-valued process such that each finite dimensional projection $(\inner{X_{t}}{\phi}: t \geq 0)$ has continuous paths, then the process has a $\Phi'$-valued version with continuous paths in $\Phi'$. Some other authors extended Mitoma's result to more general classes of nuclear spaces (see \cite{Fernique:1989, Fouque:1984, Martias:1988}) but it was only until very recently that in \cite{FonsecaMora:2018}, the author shown that Mitoma's result can be extended to a general nuclear space $\Phi$ and to a cylindrical stochastic process $(X_{t}: t \geq 0)$ that has locally equicontinuous Fourier transforms.

In this paper we carried out a further extension to the result in \cite{FonsecaMora:2018} by showing that if $(X_{t}: t \geq 0)$ is a cylindrical process in the dual $\Phi'$ of a general locally convex space such that each $(X_{t}(\phi): t \geq 0)$ has continuous (respectively c\`{a}dl\`{a}g) paths in $\R$, then local equicontinuity of its Fourier transforms in the Sazonov topology on $\Phi$ is sufficient to have a $\Phi'$-valued version $(Y_{t}: t \geq 0)$ with continuous (respectively c\`{a}dl\`{a}g) paths in $\Phi'$. We call this result again as the \emph{regularization theorem} in honour to the original work of It\^{o} and Nawata. Some other results of insterest are proven and in particular we will show that our result includes the well-known results of radonification of cylindrical processes by a Hilbert-Schmidt operator on Hilbert spaces (see \cite{BadrikianUstunel:1996, JakubowskiEtAl:2002}). As an application of our results we provide sufficient conditions for a cylindrical L\'{e}vy process in $\Phi'$ to have a regularized version that is also a $\Phi'$-valued L\'{e}vy process.  

The motivation behind our work is twofold. First, from a theoretical perspective our extended version of the regularization theorem can be though as the corresponding natural extension of the Minlos-Sazonov theorem for cylindrical measures to the context of cylindrical stochastic processes. Second,  the regularization theorem in the nuclear space setting has proven to be a very useful tool on the study of infinite dimensional stochastic processes. In particular, it has been applied to prove the L\'{e}vy-It\^{o} decomposition for L\'{e}vy processes in the dual of a nuclear space \cite{FonsecaMora:Levy}, to construct stochastic integrals driven by L\'{e}vy processes \cite{FonsecaMora:2018SPDE}, to the study of the Skokrokhod topology on the dual of a nuclear space \cite{FonsecaMora:Skorokhod} and to the study of semimartingales in duals of nuclear spaces \cite{FonsecaMora:Semi}. We hope that the results introduced in this paper can found many more applications. 

The organization of the paper is as follows. In Section \ref{sectionPrelim} we review the basic concepts on locally convex spaces, cylindrical processes and the Sazonov topology that we will need throughout this paper. Section \ref{sectionRegulTheorem} is devoted to the proof of our main result (Theorem \ref{theoRegularizationTheoremLCS}). The case of cylindrical processes with moments is studied in Section \ref{sectRegFiniMomen}. In Section \ref{sectReguThroHilbSchm} we prove that a cylindrical process can be regularized when mapped through a ``generalized'' Hilbert-Schmidt operator. In Section \ref{sectCyliLevy}
we apply our results to the study of regularied versions to cylindrical L\'{e}vy processes. Finally, in Section \ref{sectDiscuComp} we compare our results with those available on the literature.

\section{Preliminaries}\label{sectionPrelim}

\subsection{Locally Convex Spaces and Seminorms}

Let $\Phi$ be a locally convex space (we will only consider vector spaces over $\R$). 
If $p$ is a continuous semi-norm on $\Phi$ and $r>0$, the closed ball of radius $r$ of $p$ given by $B_{p}(r) = \left\{ \phi \in \Phi: p(\phi) \leq r \right\}$ is a closed, convex, balanced neighborhood of zero in $\Phi$. A continuous semi-norm (respectively a norm) $p$ on $\Phi$ is called \emph{Hilbertian} if $p(\phi)^{2}=Q(\phi,\phi)$, for all $\phi \in \Phi$, where $Q$ is a symmetric, non-negative bilinear form (respectively inner product) on $\Phi \times \Phi$. Let $\Phi_{p}$ be the Banach space (Hilbert if $p$ is Hilbertian)  that corresponds to the completion of the normed space $(\Phi / \mbox{ker}(p), \tilde{p})$, where $\tilde{p}(\phi+\mbox{ker}(p))=p(\phi)$ for each $\phi \in \Phi$. A continuous seminorm $p$ on $\Phi$ is called \emph{separable} if the space $\Phi_{p}$ is separable. A topology $\tau$ on $\Phi$ is called \emph{multi-Hilbertian} if $\tau$ is generated by a family of separable Hilbertian seminorms. 

If $p$ is a continuous seminorm on $\Phi$, the quotient map $\Phi \rightarrow \Phi / \mbox{ker}(p)$ has an unique continuous linear extension $i_{p}:\Phi \rightarrow \Phi_{p}$.  Let $q$ be another continuous semi-norm on $\Phi$ for which $p \leq q$. In this case, $\mbox{ker}(q) \subseteq \mbox{ker}(p)$. Moreover, the inclusion map from $\Phi / \mbox{ker}(q)$ into $\Phi / \mbox{ker}(p)$ is linear and continuous, and therefore it has a unique continuous extension $i_{p,q}:\Phi_{q} \rightarrow \Phi_{p}$. Furthermore, we have the following relation: $i_{p}=i_{p,q} \circ i_{q}$. The space of continuous linear operators from $\Phi_{q}$ into $\Phi_{p}$ is denoted by $\mathcal{L}(\Phi_{q},\Phi_{p})$. We employ an analogous notation for operators between the dual spaces $\Phi'_{p}$ and $\Phi'_{q}$. 

We denote by $\Phi'$ the topological dual of $\Phi$ and by $\inner{f}{\phi}$ the canonical pairing of elements $f \in \Phi'$, $\phi \in \Phi$. Unless otherwise specified, $\Phi'$ will always be consider equipped with its \emph{strong topology}, i.e. the topology on $\Phi'$ generated by the family of semi-norms $( \eta_{B} )$, where for each $B \subseteq \Phi$ bounded we have $\eta_{B}(f)=\sup \{ \abs{\inner{f}{\phi}}: \phi \in B \}$ for all $f \in \Phi'$.  If $p$ is a continuous Hilbertian semi-norm on $\Phi$, then we denote by $\Phi'_{p}$ the Hilbert space dual to $\Phi_{p}$. The dual norm $p'$ on $\Phi'_{p}$ is given by $p'(f)=\sup \{ \abs{\inner{f}{\phi}}:  \phi \in B_{p}(1) \}$ for all $ f \in \Phi'_{p}$. Moreover, the dual operator $i_{p}'$ corresponds to the canonical inclusion from $\Phi'_{p}$ into $\Phi'$ and it is linear and continuous. 

Let $( p_{n} : n \in \N)$ be an increasing sequence of separable continuous Hilbertian semi-norms on $(\Phi,\tau)$. We denote by $\theta$ the locally convex topology on $\Phi$ generated by the family $( p_{n} : n \in \N)$. The topology $\theta$ is weaker than $\tau$. We   call $\theta$ a (weaker) \emph{countably Hilbertian topology} on $\Phi$ and we denote by $\Phi_{\theta}$ the space $(\Phi,\theta)$ and by $\widetilde{\Phi_{\theta}}$ its completion. The space $\widetilde{\Phi_{\theta}}$ is a (not necessarily Hausdorff) separable, complete, pseudo-metrizable (hence Baire and ultrabornological; see Example 13.2.8(b) and Theorem 13.2.12 in \cite{NariciBeckenstein}) locally convex space and its dual space satisfies $(\widetilde{\Phi_{\theta}})'=(\Phi_{\theta})'=\bigcup_{n \in \N} \Phi'_{p_{n}}$ (see \cite{FonsecaMora:2018}, Proposition 2.4). The space $(\widetilde{\Phi_{\theta}})'$ being the inductive limit of the Souslin spaces $\Phi'_{p_{n}}$ is again a Souslin space (see \cite{Treves}, Proposition A.4(c), p.551). 

A locally convex space is called \emph{ultrabornological} if it is the inductive limit of a family of Banach spaces. For equivalent definitions see \cite{Jarchow, NariciBeckenstein}. 
  
Let us recall that a (Hausdorff) locally convex space $(\Phi,\tau)$ is called \emph{nuclear} if its topology $\tau$ is generated by a family $\Pi$ of Hilbertian semi-norms such that for each $p \in \Pi$ there exists $q \in \Pi$, satisfying $p \leq q$ and the canonical inclusion $i_{p,q}: \Phi_{q} \rightarrow \Phi_{p}$ is Hilbert-Schmidt. Other equivalent definitions of nuclear spaces can be found in \cite{Pietsch, Treves}.


\subsection{Cylindrical and Stochastic Processes} \label{subSectionCylAndStocProcess}

Let $E$ be a topological space and denote by $\mathcal{B}(E)$ its Borel $\sigma$-algebra. Recall that a Borel measure $\mu$ on $E$ is called a \emph{Radon measure} if for every $\Gamma \in \mathcal{B}(E)$ and $\epsilon >0$, there exist a compact set $K \subseteq \Gamma$ such that $\mu(\Gamma \backslash K) < \epsilon$. In general not every Borel measure on $E$ is Radon.
We denote by  $\goth{M}_{R}^{1}(E)$ the space of all Radon probability measures on $E$.  A subset $M \subseteq \goth{M}_{R}^{1}(E)$ is called \emph{uniformly tight} if for every $\epsilon >0$ there exist a compact set $K \subseteq E$ such that $\mu (K^{c})< \epsilon$ for all $\mu \in M$. A sequence $(\mu_{n}:n \in \N)$ \emph{converges weakly} to $\mu$  in $\goth{M}_{R}^{1}(E)$ if $\int_{E} f d\mu_{n} \rightarrow \int_{E} f d\mu$ for every continuous bounded function $f$ on $E$.  

Let $\Phi$ be a  locally convex space. Given $M \subseteq \Phi$, the cylindrical algebra on $\Phi'$ based on $M$ is the collection $\mathcal{Z}(\Phi',M)$ of all the \emph{cylindrical sets} of the form 
$$\mathcal{Z}\left(\phi_{1}, \dots, \phi_{n}; A \right) = \pi_{\phi_{1}, \dots, \phi_{n}}^{-1} (A),$$ 
for $n \in \N$, $\phi_{1}, \dots, \phi_{n} \in M$, $A \in \mathcal{B}\left(\R^{n}\right)$ and where $\pi_{\phi_{1}, \dots, \phi_{n}}(f)= \left(\inner{f}{\phi_{1}}, \dots, \inner{f}{\phi_{n}}\right)$. 
The $\sigma$-algebra generated by $\mathcal{Z}(\Phi',M)$ is denoted by $\mathcal{C}(\Phi',M)$. If  $M$ is finite we have $\mathcal{C}(\Phi',M)=\mathcal{Z}(\Phi',M)$. Moreover, we always have $\mathcal{C}(\Phi')\defeq \mathcal{C}(\Phi',\Phi) \subseteq \mathcal{B}(\Phi')$, but equality is not true in general.
A function $\mu: \mathcal{Z}(\Phi',\Phi) \rightarrow [0,\infty]$ is called a \emph{cylindrical measure} on $\Phi'$ if for each finite subset $M \subseteq \Phi$ the restriction of $\mu$ to $\mathcal{C}(\Phi',M)$ is a measure. A cylindrical measure $\mu$ is said to be \emph{finite} if $\mu(\Phi')< \infty$ and a \emph{cylindrical probability measure} if $\mu(\Phi')=1$. The \emph{Fourier transform}\index{characteristic function} of $\mu$ is the function $\widehat{\mu}: \Phi \rightarrow \C$ defined by 
$$ \widehat{\mu}(\phi)= \int_{\Phi'} e^{i \inner{f}{\phi}} \mu(df)= \int_{-\infty}^{\infty} \, e^{iz} \mu \circ \pi_{\phi}^{-1}(dz), \quad \forall \, \phi \in \Phi. $$

Let $\ProbSpace$ be a (complete) probability space. We denote by $L^{0} \ProbSpace$ the space of equivalence classes of real-valued random variables  defined on $\ProbSpace$.  We always consider the space $L^{0} \ProbSpace$ equipped with the topology of convergence in probability and in this case it is a complete, metrizable, topological vector space (see \cite{Badrikian}). 

A \emph{cylindrical random variable}\index{cylindrical random variable} in $\Phi'$ is a linear map $X: \Phi \rightarrow L^{0} \ProbSpace$ (see \cite{FonsecaMora:2018}).
If $Z=\mathcal{Z}\left(\phi_{1}, \dots, \phi_{n}; A \right)$ is a cylindrical set, for $\phi_{1}, \dots, \phi_{n} \in \Phi$ and $A \in \mathcal{B}\left(\R^{n}\right)$, let 
\begin{equation*} 
\mu_{X}(Z) \defeq \Prob \left( ( X(\phi_{1}), \dots, X(\phi_{n})) \in A  \right).
\end{equation*}
The map $\mu_{X}$ is a cylindrical probability measure on $\Phi'$ and it is called the \emph{cylindrical distribution} of $X$. The \emph{Fourier transform} of $X$ is defined to be the Fourier transform $\widehat{\mu}_{X}: \Phi \rightarrow \C$ of its cylindrical distribution $\mu_{X}$.


Let $X$ be a $\Phi'$-valued random variable, i.e. $X:\Omega \rightarrow \Phi'$ is a $\mathscr{F}/\mathcal{B}(\Phi')$-measurable map. We denote by $\mu_{X}$ the probability distribution of $X$, i.e. $\mu_{X}(\Gamma)=\Prob \left( X \in  \Gamma \right)$, $\forall \, \Gamma \in \mathcal{B}(\Phi')$; it is a Borel probability measure on $\Phi'$. For each $\phi \in \Phi$ we denote by $\inner{X}{\phi}$ the real-valued random variable defined by $\inner{X}{\phi}(\omega) \defeq \inner{X(\omega)}{\phi}$, for all $\omega \in \Omega$. The linear mapping $\phi \mapsto \inner{X}{\phi}$ is called the \emph{cylindrical random variable} induced/defined by $X$. 
 
Let $J=\R_{+} \defeq [0,\infty)$ or $J=[0,T]$ for  $T>0$. We say that $X=( X_{t}: t \in J)$ is a \emph{cylindrical process} in $\Phi'$ if $X_{t}$ is a cylindrical random variable for each $t \in J$. We say that $X$ is  \emph{$n$-integrable} if for every $\phi \in \Phi$, $X(\phi)=( X_{t}(\phi): t \in J)$ is $n$-integrable. 
Clearly, any $\Phi'$-valued stochastic processes $X=( X_{t}: t \in J)$ induces/defines a cylindrical process under the prescription: $\inner{X}{\phi}=( \inner{X_{t}}{\phi}: t \in J)$, for each $\phi \in \Phi$. 

If $X$ is a cylindrical random variable in $\Phi'$, a $\Phi'$-valued random variable $Y$ is a called a \emph{version} of $X$ if for every $\phi \in \Phi$, $X(\phi)=\inner{Y}{\phi}$ $\Prob$-a.e. A $\Phi'$-valued process $Y=(Y_{t}:t \in J)$ is said to be a $\Phi'$-valued \emph{version} of the cylindrical process $X=(X_{t}: t \in J)$ on $\Phi'$ if for each $t \in J$, $Y_{t}$ is a $\Phi'$-valued version of $X_{t}$.  

Let  $X=( X_{t}: t \in J)$ be a $\Phi'$-valued process.  We say that $X$ is \emph{continuous} (respectively \emph{c\`{a}dl\`{a}g}) if for $\Prob$-a.e. $\omega \in \Omega$, the \emph{sample paths} $t \mapsto X_{t}(w) \in \Phi'$ of $X$ are continuous (respectively right-continuous with left limits).

A $\Phi'$-valued random variable $X$ is called \emph{regular} if there exists a weaker countably Hilbertian topology $\theta$ on $\Phi$ such that $\Prob( \omega: X(\omega) \in (\widetilde{\Phi_{\theta}})')=1$. Furthermore, a $\Phi'$-valued process $Y=(Y_{t}:t \in J)$ is said to be \emph{regular} if $Y_{t}$ is a regular random variable $\forall t \in J$. 

Let $T>0$ and let $(E,\norm{\cdot})$ be a Banach space. We denote by $C(T,E)$ the space of $E$-valued continuous mappings on $[0,T]$. The space $C(T,E)$ is a Banach space when equipped with the topology  of uniform convergence on $[0,T]$, i.e. with the norm 
$$ F \mapsto \sup_{t \in [0,T]} \norm{ F(t)}. $$
Similarly, we denote by $C_{\Prob}(T,E)$ the space of continuous $E$-valued processes defined on $[0,T]$. This space is a complete, metrizable, topological vector space when equipped with the topology of convergence in probability uniformly on $[0,T]$. Moreover, for a given $k \in \N$, we denote by $C_{\Prob}^{k}(T,E)$ the space of continuous $E$-valued processes defined on $[0,T]$ and which have finite moment of order $k$ uniformly on $[0,T]$. The space $C_{\Prob}^{k}(T,E)$ is Banach when equipped with the norm 
$$ X \mapsto \left[ \Exp \left( \sup_{t \in [0,T]} \norm{ X_{t}}^{k} \right) \right]^{1/k}.$$
If we replace the attribute continuous by c\`{a}dl\`{a}g we can define in a complete analogue way the spaces $D(T,E)$, $D_{\Prob}(T,E)$, $D_{\Prob}^{k}(T,E)$ which possesses the same properties described above for its continuous counterparts.

\subsection{The Sazonov's topology} \label{sectSazonovTopo}

Let $(\Phi, \tau)$ be a locally convex space. In this section we recall the definition of the Sazonov topology on $\Phi$ and some of its properties. 
For further details see \cite{BogachevSmolyanovTVS, SchwartzRM, SmolyanovFomin}.

Let $\mathscr{P}(\Phi, \tau)$ denote the collection of all the seminorms on $(\Phi, \tau)$ defined in the following way: $p \in  \mathscr{P}(\Phi, \tau)$ if and only if $p$ is a continuous Hilbertian seminorm on $\Phi$ for which there exists a continuous separable Hilbertian seminorm $q$ on $\Phi$ such that $p \leq q$, and the canonical inclusion $i_{p,q}: \Phi_{q} \rightarrow \Phi_{p}$ is Hilbert-Schmidt (observe that $p$ is separable since $i_{p,q}$ is sujective and being a compact operator $i_{p,q}$ has a separable image). 

The collection $\mathscr{P}(\Phi, \tau)$ is non-empty as every seminorm on $\Phi$ continuous with respect to the weak topology  is a member of $\mathscr{P}(\Phi, \tau)$.  
The locally convex topology on $\Phi$ generated by the family of seminorms $\mathscr{P}(\Phi, \tau)$ is called the \emph{Sazonov topology} or the \emph{Hilbert-Schmidt topology} on $\Phi$ with respect to the topology $\tau$ and is denote by $\tau_{HS}$. If $\Phi$ is a separable Hilbert space with norm $\norm{\cdot}$, the Sazonov topology can be generated by the collection of all the seminorms on $\Phi$ of the form $p_{S}(\phi)=\norm{S\phi}$ $\forall \, \phi \in \Phi$, where $S$ runs over the totally of all Hilbert-Schmidt operators on $\Phi$ (see \cite{SmolyanovFomin}).

Let $\sigma$ denotes the weak topology on $(\Phi, \tau)$. Considering finite dimensional subspaces on $\Phi$ as Hilbert spaces, it is clear  that $\sigma$ is weaker than $\tau_{HS}$. On the other hand, each $p \in  \mathscr{P}(\Phi, \tau)$ is a continuous Hilbertian seminorm on $\Phi$ and therefore we have that $\tau_{HS}$ is weaker than $\tau$. The equality $\tau_{HS} = \tau$ holds if and only if $(\Phi, \tau)$ is a nuclear space. Moreover, in general $(\Phi, \tau_{HS})$ is not a nuclear space. For a counterexample see \cite{Yamasaki}, Example 18.1. 

Let $\Phi$, $\Psi$ denote two locally convex spaces. We will need the following extension of the definition of Hilbert-Schmidt operators introduced in \cite{Badrikian}. A linear operator $S$ from $\Phi$ into $\Psi$ is called a \emph{Hilbert-Schmidt operator} if there exists a continuous Hilbertian seminorm $p$ on $\Phi$ and a bounded, convex, balanced subset $B$ of $\Psi$ such that $\Psi_{B} \defeq \bigcup_{n \in \N} n B$ equipped with the norm $p_{B}(\psi)=\inf \{  \lambda >0: \psi \in \lambda B \} $ is a Hilbert space, $S(B_{p}(1))\subseteq B$, and the map $S_{0}$ in $\mathcal{L}(\Phi_{p}, \Psi_{B})$ induced by $S$ is a Hilbert-Schmidt operator. The above description means that $S$ can be decomposed as $S=j_{B} \circ S_{0} \circ i_{p}$ for $j_{B}:\Psi_{B} \rightarrow \Psi$ the canonical inclusion, i.e.  the following diagram is satisfied:
$$
\begin{gathered}
\xymatrix{ 
\Phi \ar[r]^{S} \ar[d]_{i_{p}} & \Psi  \\
\Phi_{p} \ar[r]_{S_{0}} & \Psi_{B} \ar[u]_{j_{B}}
}
\end{gathered}
$$
If $\Psi$ is a Hilbert space, then $\Psi_{B} = \Psi$ and  
$j_{B}$ is an isometry on $\Psi$. If furthermore, $\Phi$ is also a  Hilbert space, then $\Phi_{p} = \Phi$ and $i_{p}$ is an isometry on $\Phi$, hence $S$ becomes a Hilbert-Schmidt operator in the usual sense.  

It is proven in \cite{Badrikian} (Expos\'{e} no.10, Proposition 4, p.168) that if $\Phi$ is a locally convex space and $\Psi$ is a Hilbert space, $R \in \mathcal{L}(\Phi,\Psi)$ is Hilbert-Schmidt if and only if is continuous when $\Phi$ is equipped in its Sazonov topology $\tau_{HS}$. Hence, if $\Psi$ is locally convex and $S \in \mathcal{L}(\Phi,\Psi)$ is Hilbert-Schmidt it is continuous from $(\Phi,\tau_{HS})$ into $\Psi$.

\section{The Regularization Theorem for Cylindrical Stochastic Processes}\label{sectionRegulTheorem}

\begin{assu}
From now on and unless otherwise specified, $\Phi$ will always denote a locally convex space. 
\end{assu}

Our main result is the following that provides sufficient conditions for a cylindrical process in $\Phi'$ to have a $\Phi'$-valued continuous or c\`{a}dl\`{a}g version. 

\begin{theo}[Regularization Theorem in a Locally Convex Space]\label{theoRegularizationTheoremLCS}
Let $X=(X_{t}: t \geq 0)$ be a cylindrical process in $\Phi'$ satisfying:
\begin{enumerate}
\item For each $\phi \in \Phi$, the real-valued process $X(\phi)=( X_{t}(\phi) : t \geq 0)$ has a continuous (respectively c\`{a}dl\`{a}g) version.
\item For every $T > 0$, the Fourier transforms 
of the family $( X_{t}: t \in [0,T] )$ are equicontinuous (at the origin) in the Sazonov topology $\tau_{HS}$.
\end{enumerate}
Then, there exists a countably Hilbertian topology $\vartheta$ on $\Phi$ and a $(\widetilde{\Phi_{\vartheta}})'$-valued continuous (respectively c\`{a}dl\`{a}g) process $Y= (Y_{t}: t \geq 0)$, such that for every $\phi \in \Phi$, $\inner{Y}{\phi}= ( \inner{Y_{t}}{\phi} : t \geq 0)$ is a version of $X(\phi)= ( X_{t}(\phi) : t \geq 0)$. Moreover, $Y$ is a $\Phi'$-valued, regular, continuous (respectively c\`{a}dl\`{a}g) version of $X$ that is unique up to indistinguishable versions. 
\end{theo}

The proof of Theorem \ref{theoRegularizationTheoremLCS} is a modification of the arguments in the proof of Theorem 3.2 in \cite{FonsecaMora:2018}. Hence, we sometimes refer the reader to \cite{FonsecaMora:2018} when there will be arguments that follows similarly to those used there. 

First, to prove the uniquess up to indistinguishable versions, suppose that $Y$ and $Z$ are two $\Phi'$-valued processes satisfying the conclusion of Theorem \ref{theoRegularizationTheoremLCS}. Then for every $\phi \in \Phi$ and $t \geq 0$, $\inner{Y_{t}}{\phi}=X_{t}(\phi)=\inner{Z_{t}}{\phi}$ $\Prob$-a.e., and because $Y$ and $Z$ are both  $\Phi'$-valued, regular, continuous (respectively c\`{a}dl\`{a}g)processes it follows from Proposition 2.12 in  \cite{FonsecaMora:2018} that $Y$ and $Z$ are indistinguishable processes. 

Now, for the proof of Theorem \ref{theoRegularizationTheoremLCS}, it is enough to show that the result is valid for a cylindrical processed defined on the bounded interval of time $[0,T]$, for $T>0$, and under the assumption that each  $X(\phi)= ( X_{t}(\phi) : t \in [0,T])$ has a continuous  version. The arguments are similar to those used in the proof of Theorem 3.2 in \cite{FonsecaMora:2018}, but for the reader's convenience we summarize the main steps.

In effect, if the result is valid for every $T>0$, then for every $n \in \N$ we can find a weaker countably Hilbertian topology $\vartheta_{n}$ on $\Phi$ and a  $(\widetilde{\Phi_{\vartheta_{n}}})'$-valued continuous (respectively c\`{a}dl\`{a}g) process $Y^{(n)}= \left(Y^{(n)}_{t}: t \in [0,n] \right)$ such that for each $\phi \in \Phi$, $\inner{Y^{(n)}}{\phi}= \left( \inner{Y^{(n)}_{t}}{\phi} : t \in [0,n] \right)$ is a version of $X(\phi)= ( X_{t}(\phi) : t \in [0,n])$.

Let $\vartheta$ denote the countably Hilbertian topology on $\Phi$ generated by the families of seminorms generating the topologies $\vartheta_{n}$, $n \in \N$. The topology $\vartheta$ is finner than each $\vartheta_{n}$, but is weaker than the given topology on $\Phi$. Therefore, through the canonical inclusion from  $(\widetilde{\Phi_{\vartheta_{n}}})'$ into  $(\widetilde{\Phi_{\vartheta}})'$, each $Y^{(n)}$ can be considered as a $(\widetilde{\Phi_{\vartheta}})'$-valued continuous process.  Since for each $n \in \N$ and $\phi \in \Phi$, the process $\inner{Y^{(n)}}{\phi}$ is a version of $X(\phi)$ on the time interval $[0,n]$, then for each $\phi \in \Phi$, $\inner{Y^{(n)}}{\phi}$ and $\inner{Y^{(n+1)}}{\phi}$ are indistinguishable as processes defined on $[0,n]$.  Then, the fact that as a $\Phi'$-valued process, each $Y^{(n)}$ is a regular process with continuous trajectories implies that $Y^{(n)}$ and $Y^{(n+1)}$ are indistinguishable as processes defined on $[0,n]$ (see Proposition 2.12 in \cite{FonsecaMora:2018}). 
 
Take $Y=\{ Y_{t} \}_{t \geq 0}$ defined by the prescription $Y_{t}=Y^{(n)}_{t}$ if $t \in [0,n]$. From the arguments in the above paragraph it is clear that $Y$ is a $(\widetilde{\Phi_{\vartheta}})'$-valued process  with continuous trajectories  such that for every $\phi \in \Phi$, $\inner{Y}{\phi}= ( \inner{Y_{t}}{\phi} : t \geq 0)$ is a version of $X(\phi)= ( X_{t}(\phi) : t \geq 0)$. The same results are obtained in the c\`{a}dl\`{a}g version case. 

From now on we will fix $T>0$ and show that the conclusions of Theorem \ref{theoRegularizationTheoremLCS} are valid for a cylindrical process defined on $[0,T]$  under the continuous version assumption.  

For each $\phi$, denote by $\widehat{X}(\phi)=( \widehat{X}_{t}(\phi) : t \in [0,T])$  a continuous version of $X(\phi)=( X_{t}(\phi) : t \in [0,T])$. Clearly, $\widehat{X}$ determines a cylindrical process $\widehat{X}=( \widehat{X}_{t}: t \in [0,T] )$ in $\Phi'$ and its Fourier transforms are equicontinuous in the Sazonov topology $\tau_{HS}$ on $\Phi$. Similar arguments to those used in the proof of Lemma 3.4 in \cite{FonsecaMora:2018} shows that there exists a countably Hilbertian topology $\theta$ on $\Phi$, generated by a sequence of  Hilbertian seminorms continuous on the Sazonov topology $\tau_{HS}$, such that the family $( \widehat{X}_{t}: t \in [0,T] )$ of mappings from $\Phi$ into $L^{0} \ProbSpace$ is $\theta$-equicontinuous.
Then, for each $t \in [0,T]$ the map $\widehat{X}_{t}$ has a continuous linear extension from $\widetilde{\Phi_{\theta}}$ into $L^{0} \ProbSpace$ (see \cite{Jarchow}, Theorem 2, Section 3.4, p.61-2). We denote this extension again by  $\widehat{X}_{t}$. Then, by following a closed graph theorem argument similar to that in the proof of Lemma 3.5 in  \cite{FonsecaMora:2018}, one can show that the linear map $\widehat{X}$ from $\widetilde{\Phi_{\theta}}$ into $C_{\Prob}(T,\R)$ given by $\phi \mapsto  \widehat{X}(\phi)=( \widehat{X}_{t}(\phi) : t \in [0,T])$ is continuous.  
Then, as in Lemma 3.7 in \cite{FonsecaMora:2018} we get the following result:

\begin{lemm}\label{lemmPnContiCharacFunctSupremum} For every $\epsilon >0$ there exists a  $\theta$-continuous (hence $\tau_{HS}$-continuous) Hilbertian semi-norm $p$ on $\Phi$ such that 
\begin{equation} \label{pnContinuityCharactFuncSupremum}
\Exp \left( \sup_{t \in [0,T]} \abs{1-e^{i \widehat{X}_{t}(\phi)}} \right) \leq \epsilon + 2p(\phi)^{2}, \quad \forall \, \phi \in \Phi. 
\end{equation}
\end{lemm}

Let $( \epsilon_{n}: n \in \N)$ be a sequence of positive numbers such that $\lim_{n \rightarrow \infty} \epsilon_{n}=0$. From Lemma \ref{lemmPnContiCharacFunctSupremum} there exists an increasing sequence of $\theta$-continuous (hence $\tau_{HS}$-continuous) Hilbertian semi-norms $( p_{n} : n \in \N)$ on $\Phi$ such that  $\epsilon_{n}$ and $p_{n}$ satisfy  \eqref{pnContinuityCharactFuncSupremum} $\forall n \in \N$.

Since each $p_{n}$ is $\tau_{HS}$-continuous, there exists a sequence (which we can choose increasing) of separable continuous Hilbertian seminorms $(q_{n}: n \in \N)$ on $\Phi$ such that $\forall n \in \N$, $p_{n} \leq q_{n}$, and the inclusion $i_{p_{n},q_{n}}: \Phi_{q_{n}} \rightarrow \Phi_{p_{n}}$ is Hilbert-Schmidt. 

Denote by $\vartheta$ the countably Hilbertian topology on $\Phi$ generated by the seminorms $(q_{n}: n \in \N)$. The space $\Phi_{\vartheta}$ is separable.  Let $( \xi_{k}: k \in \N )$  be a countable dense subset of $\Phi_{\vartheta}$. For every $n \in \N$, from an application of the Schmidt orthogonalization procedure we can find a complete orthonormal system $( \phi_{j}^{q_{n}} : j \in \N) \subseteq \Phi$ of $\Phi_{q_{n}}$, such that 
\begin{equation} \label{decompDenseSetInTermsOrtoBasis}
\xi_{k}= \sum_{j=1}^{k} a_{j,k,n} \, \phi_{j}^{q_{n}} + \varphi_{k,n}, \quad \forall \, k \in \N, 
\end{equation}
with $a_{j,k,n} \in \R$ and $\varphi_{k,n} \in \mbox{Ker}(q_{n})$, for each $j,k \in \N$. 

Let $n \in \N$. From \eqref{pnContinuityCharactFuncSupremum} and by following similar calculations to those used in the proof of Lemma 3.8 in \cite{FonsecaMora:2018}, p.878-9, one can show that: 
\begin{equation} \label{firstPartInequalityProofContCadlagVersionLCS}
\Prob \left( \sup_{t \in [0,T]} \sum_{j=1}^{\infty} \abs{\widehat{X}_{t}(\phi_{j}^{q_{n}})}^{2} <  \infty \right)  \geq 1- \frac{\sqrt{e}}{\sqrt{e}-1} \epsilon_{n}. 
\end{equation} 
\begin{equation} \label{secondPartInequalityProofContCadlagVersionLCS}
\Prob \left( \sup_{t \in [0,T]} \sum_{j=1}^{\infty} \abs{\widehat{X}_{t}(\varphi_{j,n})}^{2} >0 \right)  \leq \frac{\sqrt{e}}{\sqrt{e}-1} \epsilon_{n}. 
\end{equation}
Then, if we define
\begin{equation} \label{defiSetBNRegulaTheo}
B_{n}=\left\{ \omega: \sup_{t \in [0,T]} \sum_{j=1}^{\infty} \abs{\widehat{X}_{t}(\phi_{j}^{q_{n}})(\omega)}^{2} <  \infty \mbox{ and } \widehat{X}_{t}(\varphi_{j,n})(\omega) =0,  \forall  t \in [0,T],  j \in \N \right\} 
\end{equation}
it follows from \eqref{firstPartInequalityProofContCadlagVersionLCS} and \eqref{secondPartInequalityProofContCadlagVersionLCS} that 
\begin{equation*} \label{probSetOmegaNRegulaTheo} 
\Prob \left( B_{n} \right)  \geq 1- 2 \frac{\sqrt{e}}{\sqrt{e}-1} \epsilon_{n}. 
\end{equation*}
For each $n \in \N$, define the sets 
\begin{equation} \label{defiSetGammaNRegulTheo}
\Gamma_{n} =  \left\{ \omega: \forall \, j \in \N, \, t \mapsto \widehat{X}_{t}(\phi_{j}^{q_{n}})(\omega) \mbox{ is continuous} \right\}. 
\end{equation} 
\begin{equation} \label{defiSetANRegulTheo}
A_{n} = \left\{ \omega: \widehat{X}_{t}(\xi_{k})(\omega)= \sum_{j=1}^{k} a_{j,k,n} \widehat{X}_{t}(\phi_{j}^{q_{n}})(\omega) + \widehat{X}_{t}(\varphi_{k,n})(\omega), \forall k \in \N, t \in [0,T]  \right\}
\end{equation}
The fact that $\widehat{X}$ is linear and continuous from $\widetilde{\Phi_{\theta}}$ into $C_{\Prob}(T,\R)$, together with \eqref{decompDenseSetInTermsOrtoBasis} shows that 
$\Prob(\Gamma_{n} \cap A_{n})=1$ for every $n \in \N$. Therefore, if we define 
\begin{equation} \label{defiSetLambdaNRegulTheo}
\Lambda_{n} = B_{n} \cap \Gamma_{n} \cap A_{n}. 
\end{equation}
we have 
\begin{equation} \label{probSetLambdaNRegulTheo}
\Prob \left( \Lambda_{n} \right) \geq 1- 2 \frac{\sqrt{e}}{\sqrt{e}-1} \epsilon_{n}. 
\end{equation}
Since $\lim_{n \rightarrow \infty} \epsilon_{n}=0$, it follows from \eqref{probSetLambdaNRegulTheo} that $\Prob \left(\bigcup_{n} \Lambda_{n} \right) =1$. 
Set $\Omega_{1}=\Lambda_{1}$, and $\Omega_{n}=\Lambda_{n} \setminus \Lambda_{n-1}$ for $n \geq 2$. Then, $\Prob \left(\bigcup_{n} \Omega_{n} \right) =1$.  

For each $n \in \N$, let $( f_{j}^{q_{n}}: j \in \N)$ be a complete orthonormal system in $\Phi'_{q_{n}}$ dual to $(\phi_{j}^{q_{n}} : j \in \N)$, i.e. $\inner{f_{j}^{q_{n}}}{\phi_{i}^{q_{n}}}=\delta_{i,j}$ where $\delta_{i,j}=1$ if $i=j$ and $\delta_{i,j}=0$ if $i \neq j$. 
For each $t \in [0,T]$, we define 
\begin{equation} \label{defiVersionYRegularTheorem}
Y_{t}(\omega) \defeq 
\begin{cases}
\sum_{j=1}^{\infty} \widehat{X}_{t}(\phi_{j}^{q_{n}})(\omega) f_{j}^{q_{n}}, & \mbox{if } \omega \in \Omega_{n}, \\
0, & \mbox{if } \omega \notin \bigcup_{n} \Omega_{n}. 
\end{cases}
\end{equation}
We will show that $Y= (Y_{t}: t \geq 0)$ is a well-defined $(\widetilde{\Phi_{\vartheta}})'$-valued process with continuous trajectories. To do this, suppose $\omega \in \Omega_{n}$. Since $\omega \in \Gamma_{n}$, then for every $m \in \N$,  
$$ \sum_{j=1}^{m} \widehat{X}_{t}(\phi_{j}^{q_{n}})(\omega) f_{j}^{q_{n}} \in C(T,\R) \otimes \Phi'_{q_{n}} \subseteq C(T,\Phi'_{q_{n}}).$$
Moreover, for $m \geq k \geq 1$, it follows from \eqref{defiSetBNRegulaTheo} that 
\begin{flalign*} 
& \sup_{t \in [0,T]} q_{n}'\left(\sum_{j=1}^{m} \widehat{X}_{t}(\phi_{j}^{q_{n}})(\omega) f_{j}^{q_{n}}-\sum_{j=1}^{k} \widehat{X}_{t}(\phi_{j}^{q_{n}})(\omega) f_{j}^{q_{n}} \right)^{2} \\
& = \sup_{t \in [0,T]} \sum_{j =k+1}^{m} \abs{\widehat{X}_{t}(\phi_{j}^{q_{n}})(\omega)}^{2} \rightarrow 0, \quad \mbox{ as } m,k \rightarrow \infty. 
\end{flalign*}
Since the Banach space $ C(T,\Phi'_{q_{n}})$ is complete, 
the Cauchy sequence $\sum_{j=1}^{m} \widehat{X}_{t}(\phi_{j}^{q_{n}})(\omega) f_{j}^{q_{n}}$, $m \in \N$, converges and its limit $\sum_{j=1}^{\infty} \widehat{X}_{t}(\phi_{j}^{q_{n}})(\omega) f_{j}^{q_{n}}$ belongs to $ C(T,\Phi'_{q_{n}})$. 
Then, it follows from \eqref{defiVersionYRegularTheorem} that for $\omega \in \Omega_{n}$, $Y_{t}(\omega) \in \Phi'_{q_{n}}$ and the map $t \mapsto Y_{t}(\omega)$ is continuous in $\Phi'_{q_{n}}$. Since the canonical inclusion from $\Phi'_{q_{n}}$ into  $(\widetilde{\Phi_{\vartheta}})'$ is linear and continuous, and since $\Prob \left( \cup_{n} \Omega_{n} \right)=1$, we conclude that $Y$ is a $(\widetilde{\Phi_{\vartheta}})'$-valued process with continuous trajectories. 

Our final task is to check that for each $\phi \in \Phi$, $\inner{Y}{\phi}= ( \inner{Y_{t}}{\phi} : t \geq 0)$ is a version of $X(\phi)= ( X_{t}(\phi) : t \geq 0)$. 
Let $\phi \in \Phi$. Since $(\xi_{k}: k \in \N )$ is dense in $\Phi_{\vartheta}$, there exists a subsequence $(\xi_{k_{j}}: j \in \N )$ such that $\xi_{k_{j}} \rightarrow \phi$ in $\Phi_{\vartheta}$. 

Now, observe that for each $k \in \N$, we have $\inner{Y_{t}(\omega)}{\xi_{k}}=\widehat{X}_{t}(\xi_{k})(\omega)$ for $\omega \in \cup_{n} \Omega_{n}$. In effect, if $\omega \in \Omega_{n}$, by \eqref{defiVersionYRegularTheorem} and because $\inner{f_{j}^{q_{n}}}{\phi_{i}^{q_{n}}}=\delta_{i,j}$, it is clear that 
$\inner{Y_{t}(\omega)}{\phi_{j}^{q_{n}}}=\widehat{X}_{t}(\phi_{j}^{q_{n}})(\omega)$, for all $j \in \N$ and $t \in [0,T]$. Similarly, from the fact that $q_{n}(\varphi_{j,n})=0$ for all $j \in \N$, we have $\abs{\inner{f_{j}^{q_{n}}}{\varphi_{j,n}}} \leq q'_{n}(f_{j}^{q_{n}})q_{n}( \varphi_{j,n})=0$ for all $j \in \N$. Then, as $\omega \in B_{n} \cap \Gamma_{n}$ we obtain that  $\inner{Y_{t}(\omega)}{ \varphi_{j,n}}=\widehat{X}_{t}(\varphi_{j,n})(\omega)=0$ for all $t \in [0,T]$, $j \in \N$. Therefore, from \eqref{decompDenseSetInTermsOrtoBasis},   \eqref{defiSetANRegulTheo}, \eqref{defiVersionYRegularTheorem}, we have 
\begin{equation} \label{equalityInDenseSetXAndVersionYn}
\forall \, \omega \in \Omega_{n}, \quad \inner{Y_{t}(\omega)}{ \xi_{k}}=\widehat{X}_{t}(\xi_{k})(\omega), \quad \forall \, k \in \N, \, t \in [0,T]. 
\end{equation}  
On the other hand, since the map  $\widehat{X}$ from $\widetilde{\Phi_{\theta}}$ into $C_{\Prob}(T,\R)$ is linear and continuous, and since the canonical inclusion from 
$\widetilde{\Phi_{\vartheta}}$ into $\widetilde{\Phi_{\theta}}$ is linear continuous (this because the topology $\vartheta$ is finner than $\theta$), then $\widehat{X}(\xi_{k}) \rightarrow \widehat{X}(\phi)$ in  $C_{\Prob}(T,\R)$. Hence, there exists 
$\Delta_{\phi} \subseteq \Omega$ with $\Prob (\Delta_{\phi})=1$ and a subsequence $( \xi_{k_{j,\nu}}: \nu \in \N)$ of $(\xi_{k_{j}}: j \in \N)$ such that for all $\omega \in \Delta_{\phi}$, $\widehat{X}_{t}(\xi_{k_{j,\nu}})(\omega) \rightarrow \widehat{X}_{t}(\phi)(\omega)$, as $\nu \rightarrow \infty$, for all $t \in [0,T]$.

Then, \eqref{equalityInDenseSetXAndVersionYn} and the uniqueness of limits implies that 
\begin{equation} \label{equalityPaeXAndVersionYnOnPhi}
\forall \, \omega \in \Omega_{n} \cap \Delta_{\phi}, \quad  \inner{Y_{t}(\omega)}{ \phi}=\widehat{X}_{t}(\phi)(\omega), \quad \forall \, t \in [0,T]. 
\end{equation}  
Hence, as $\Prob \left( \bigcup_{n} \Omega_{n} \right)=1$ and $\widehat{X}(\phi)$ is a version of $X(\phi)$, we have proven that  $\inner{Y}{\phi}$ is a version of $X(\phi)$. This finalizes the proof of Theorem \ref{theoRegularizationTheoremLCS}.  

\begin{rema}
One can observe from the proof of Theorem \ref{theoRegularizationTheoremLCS} that the process $Y$ can be chosen such that for every $\omega \in \Omega$ and $T>0$ there exists a separable  continuous Hilbertian semi-norm $q=q(\omega,T)$ on $\Phi$ such that the map $t \mapsto Y_{t}(\omega)$ is continuous (respectively c\`{a}dl\`{a}g) from $[0,T]$ into the Hilbert space $\Phi'_{q}$. 
\end{rema}

\begin{coro} \label{coroRegulSingleSeminorm}
Let $X=(X_{t}: t \geq 0)$ be a cylindrical process in $\Phi'$ satisfying:
\begin{enumerate}
\item For each $\phi \in \Phi$, the real-valued process $X(\phi)=( X_{t}(\phi) : t \geq 0)$ has a continuous (respectively c\`{a}dl\`{a}g) version.
\item There exists a $\tau_{HS}$-continuous Hilbertian semi-norm $p$ on $\Phi$ such that for every $t \geq 0$, the Fourier transform of $X_{t}$ is $p$-continuous (at the origin). 
\end{enumerate}
Then, there exists a separable continuous Hilbertian semi-norm $q$ on $\Phi$, $p \leq q$, such that $i_{p,q}$ is Hilbert-Schmidt and a $\Phi'_{q}$-valued continuous (respectively c\`{a}dl\`{a}g) process $Y= ( Y_{t} : t \geq 0)$, such that for every $\phi \in \Phi$, $\inner{Y}{\phi}= ( \inner{Y_{t}}{\phi} : t \geq 0)$ is a version of $X(\phi)= (X_{t}(\phi) : t \geq 0)$. Moreover, $Y$ is unique up to indistinguishable versions in $\Phi'$. 
\end{coro}
\begin{prf}
As in the proof of Theorem \ref{theoRegularizationTheoremLCS}, it is enough to prove the result for a cylindrical process defined on $[0,T]$ for $T>0$ and in the continuous version case. Now, our assumptions on $X$ imply that for each $t \in [0,T]$, the map  $X_{t}:\Phi \rightarrow L^{0} \ProbSpace$ is $p$-continuous. This in turn implies that the map $\phi \mapsto (X_{t}(\phi) : t \in [0,T])$ has an extension $\widehat{X}$ that is linear and continuous from $\Phi_{p}$ into $C_{\Prob}(T,\R)$. Let $q$ be a separable continuous Hilbertian semi-norm $q$ on $\Phi$, $p \leq q$, such that $i_{p,q}$ is Hilbert-Schmidt. Then, one can follow exactly the same arguments as in the proof of Theorem \ref{theoRegularizationTheoremLCS} by taking $p_{n}=p$ and $q_{n}=q$, and in that case one conclude $B_{n}=B_{m}$, $\Gamma_{n}=\Gamma_{m}$, $A_{n}=A_{m}$ and $\Omega_{n}=\Omega_{m}$ for each $n, m \in \N$. Then, \eqref{defiVersionYRegularTheorem} defines a $\Phi'_{q}$-valued continuous process satisfying the conclusions in the statement of Corollary \ref{coroRegulSingleSeminorm}. 
\end{prf}

\section{Regularization For Cylindrical Processes With Finite Moments} \label{sectRegFiniMomen}

This section is devoted to the study of the regularization of cylindrical processes which posseses moments of finite order uniformly on a bounded interval of time. As the next result shows, in such a case one can get a considerable improvement on the results obtained in Theorem \ref{theoRegularizationTheoremLCS}. 

\begin{theo} \label{theoRegulFiniteMomentsLCS}
Let $X=(X_{t}: t \geq 0)$ be a cylindrical process in $\Phi'$ satisfying:
\begin{enumerate}
\item For each $\phi \in \Phi$, the real-valued process $X(\phi)=( X_{t}(\phi) : t \geq 0)$ has a continuous (respectively c\`{a}dl\`{a}g) version.
\item For every $T > 0$, the Fourier transforms 
of the family $( X_{t}: t \in [0,T] )$ are equicontinuous (at the origin) in the Sazonov topology $\tau_{HS}$. 
\item There exists $r \geq 2$ such that  $\Exp \left( \sup_{t \in [0,T]} \abs{X_{t}(\phi)}^{r} \right) < \infty$, $\forall \, T>0$, $\phi \in \Phi$. 
\end{enumerate}
Then, there exists an increasing sequence $(q_{n}:n \in \N)$ of separable continuous Hilbertian seminorms on $\Phi$, and a $\Phi'$-valued regular, continuous (respectively c\`{a}dl\`{a}g) version $Y= (Y_{t}: t \geq 0)$ of $X$ (unique up to indistinguishable versions), such that for every $T>0$, there exists $n \in \N$ such that $( Y_{t} : t \in [0,T])$ is a $\Phi'_{q_{n}}$-valued continuous (respectively c\`{a}dl\`{a}g) process and moreover we have  
$\Exp \left( \sup_{t \in [0,T]} q_{n}'(Y_{t})^{r} \right) < \infty.$  
\end{theo}
\begin{prf}
It is sufficient to prove the result for a given $T>0$ and the cylindrical process $(X_{t}: t \in [0,T])$ under the continuous version case. 

First, for each $\phi \in \Phi$, denote by $\widehat{X}(\phi)=(\widehat{X}_{t}(\phi): t \in [0,T])$  a continuous version of $X(\phi)=(X_{t}(\phi): t \in [0,T])$. Then, by following similar arguments to those used in the proof of Theorem \ref{theoRegularizationTheoremLCS} we can show that the map $\widehat{X}$ from $(\Phi, \tau_{HS})$ into the space $C_{\Prob}^{r}(T,\R)$ (see Section \ref{subSectionCylAndStocProcess}) given by $\phi \mapsto \widehat{X}(\phi)$ is linear and continuous. Then, if we define 
$\varrho: \Phi \mapsto [0,\infty)$ by 
$$ \varrho(\phi)= \left[ \Exp \left( \sup_{t \in [0,T]} \abs{\widehat{X}_{t}(\phi)}^{r} \right) \right]^{1/r}, \quad \forall \phi \in \Phi, $$
it follows that $\varrho$ is a $\tau_{HS}$-continuous seminorm on $\Phi$. Let $p$ be a $\tau_{HS}$-continuous Hilbertian seminorm on $\Phi$ such that $\varrho \leq p$. Then, $\widehat{X}$ has a continuous and linear extension from $\Phi_{p}$ into $C_{\Prob}^{r}(T,\R)$. 

Let $q$ be a separable continuous Hilbertian seminorm on $\Phi$ such that $p \leq q$, and the inclusion $i_{p,q}$ is Hilbert-Schmidt. We choose a countable dense subset  $( \xi_{k}: k \in \N ) \subseteq \Phi$ of $\Phi_{q}$, a complete orthonormal system $( \phi_{j}^{q} : j \in \N) \subseteq \Phi$ of $\Phi_{q}$, $a_{j,k} \in \R$ and $\varphi_{k} \in \mbox{Ker}(q)$, satisfying \eqref{decompDenseSetInTermsOrtoBasis} (observe that in our present case there is no dependence on $n$). 

Let $( f_{j}^{q}: j \in \N)$ be a complete orthonormal system in $\Phi'_{q}$ dual to $(\phi_{j}^{q} : j \in \N)$. For each $n \in \N$, let $Y^{n}$ defined by 
$$Y_{t}^{n}= \sum_{j=1}^{n} \widehat{X}_{t}(\phi_{j}^{q}) f_{j}^{q}, \quad \forall t \in [0,T].$$
It is clear that $Y^{n} \in C_{\Prob}^{r}(T,\Phi'_{q})$. Moreover, for $m >n \geq 1$, 
\begin{eqnarray*}
\Exp \left( \sup_{t \in [0,T]} q'(Y^{m}_{t}-Y^{n}_{t})^{r} \right) & \leq & \Exp \left( \sup_{t \in [0,T]} \sum_{j=n+1}^{m} \abs{\widehat{X}_{t}(\phi_{j}^{q})}^{r} \right)  \\
& \leq & \sum_{j=n+1}^{m} \Exp \left( \sup_{t \in [0,T]} \abs{\widehat{X}_{t}(\phi_{j}^{q})}^{r} \right) \\
& \leq & \sum_{j=n+1}^{m}  p(i_{p,q} \phi_{j}^{q})^{r}.   
\end{eqnarray*}
Now, since $i_{p,q}$ is Hilbert-Schmidt it belongs to the $r$-th Schatten-von Neumann class of operators (see \cite{DiestelJarchowTonge},  Proposition 4.5 and Corollary 4.8), and hence $\displaystyle{\sum_{j=1}^{\infty}  p(i_{p,q} \phi_{j}^{q})^{r}< \infty}$ (see \cite{DiestelJarchowTonge}, Theorem 4.7). 
The above calculations show that $(Y^{n}: n \in \N)$ is a Cauchy sequence in $C_{\Prob}^{r}(T,\Phi'_{q})$ and hence its limit $Y$ given by  
$$Y_{t} = \sum_{j=1}^{\infty} \widehat{X}_{t}(\phi_{j}^{q}) f_{j}^{q}, \quad \forall t \in [0,T],$$
is also an element of $C_{\Prob}^{r}(T,\Phi'_{q})$. The final step is to prove that for each $\phi \in \Phi$, $\inner{Y}{\phi}=(\inner{Y_{t}}{\phi}: t \in [0,T])$ is a version of $\widehat{X}(\phi)=(\widehat{X}_{t}(\phi): t \in [0,T])$. This can be done by following similar arguments to those used in the proof of Theorem \ref{theoRegularizationTheoremLCS} by first showing that $\inner{Y}{\xi_{k}}$ is a version of $\widehat{X}(\xi_{k})$ for each $k \in \N$, and then to use the denseness of $(\xi_{k}: k \in \N)$ on $\Phi_{q}$ and the fact that $\widehat{X}$ can be extended to be a continuous operator on $\Phi_{q}$ to get the desired conclusion. We leave the details to the reader. 
\end{prf}

\begin{exam} \label{examCyliMartingales}
Let $M=(M_{t}: t \geq 0)$ be a cylindrical martingale in $\Phi'$, i.e. $M$ is a cylindrical process such that for each $\phi \in \Phi$, $M(\phi)=(M_{t}(\phi): t \geq 0)$ is a real-valued martingale. Suppose that for each $t \geq 0$, the Fourier transform of $M_{t}$ is continuous  (at the origin) in the Sazonov topology $\tau_{HS}$. One can replicate the arguments used in the proof of Theorem 5.2 in \cite{FonsecaMora:2018} to show that $M$ satisfies the conditions $(1)$ and $(2)$ in Theorem \ref{theoRegularizationTheoremLCS}, and hence $M$ has a $\Phi'$-valued, regular, c\`{a}dl\`{a}g version $\widetilde{M}=(\widetilde{M}_{t}: t \geq 0)$. In general $\widetilde{M}$ is only cilindrically a martingale, i.e. the cylindrical process induced by $\widetilde{M}$ is a cylindrical martingale. However, if $\Phi$ is a separable Hilbert space one can show that $\widetilde{M}$ is indeed a $\Phi'$-valued martingale. 

Now, if the cylindrical martingale $M$ is $r$-integrable  for $r \geq 2$, then it is a consequence of Doob's inequality for real-valued martingales that $M$ satisfies the condition $(3)$ in Theorem \ref{theoRegulFiniteMomentsLCS}. Therefore, for any $T>0$, there exists a separable continuous Hilbertian seminorm  $q_{T}$ on $\Phi$ such that $(\widetilde{M}_{t}: t \in [0,T])$ is a $\Phi'_{q_{T}}$-valued c\`{a}dl\`{a}g martingale satisfying 
$\Exp \left( \sup_{t \in [0,T]} q_{T}'(\widetilde{M}_{t})^{r} \right) < \infty.$ 
\end{exam}

\section{Regularization Through a Hilbert-Schmidt Operator}\label{sectReguThroHilbSchm}

Suppose that $X=(X_{t}: t \geq 0)$ is a cylindrical process in the dual $\Psi'$ of a locally convex space $\Psi$. If $S$ is a linear operator from $\Phi$ into  $\Psi$, then it is clear that $X\circ S= ( X_{t}\circ S : t \geq 0)$ is a cylindrical process in $\Phi'$. Sufficient conditions for the existence of a $\Phi'$-valued continuous or c\`{a}dl\`{a}g version for $X\circ S$ are given in the next result. 

\begin{theo}  \label{theoRadonSingleHilSchOper}
Let $\Psi$ be a locally convex space equipped with a multi-Hilbertian topology and let $X=(X_{t}: t \geq 0)$ be a cylindrical process in $\Psi'$ satisfying:
\begin{enumerate}
\item For each $\psi \in \Psi$, the real-valued process $X(\psi)=( X_{t}(\psi) : t \geq 0)$ is continuous (respectively c\`{a}dl\`{a}g).
\item For every $T > 0$, the Fourier transforms 
of the family $( X_{t}: t \in [0,T] )$ are equicontinuous (at the origin) on $\Psi$.
\end{enumerate}
Let $S$ be a Hilbert-Schmidt operator from $\Phi$ into $\Psi$. 
Then, there exists a countably Hilbertian topology $\vartheta$ on $\Phi$ and a $(\widetilde{\Phi_{\vartheta}})'$-valued continuous (respectively c\`{a}dl\`{a}g) process $Y= (Y_{t}: t \geq 0)$, such that for every $\phi \in \Phi$, $\inner{Y}{\phi}= ( \inner{Y_{t}}{\phi} : t \geq 0)$ is a version of $X\circ S(\phi)= ( X_{t}\circ S ( \phi) : t \geq 0)$. Moreover, $Y$ is a $\Phi'$-valued, regular, continuous (respectively c\`{a}dl\`{a}g) version of $X \circ S$ that is unique up to indistinguishable versions. 
\end{theo}
\begin{prf} We will show that $X \circ S$ satisfies the conditions in Theorem \ref{theoRegularizationTheoremLCS}. We do this under the continuous version assumption. 

The condition \emph{(1)} in Theorem \ref{theoRegularizationTheoremLCS} is obviously satisfied. 
Let $T>0$. As in the proof of Theorem \ref{theoRegularizationTheoremLCS}, $X$ induces a linear map $\psi \mapsto X(\psi)=(X_{t}(\psi): t \in [0,T])$ from $\Psi$ into $C_{\Prob}(T,\R)$. Our assumption \emph{(2)} and a closed graph theorem argument similar to that used in the proof of Lemma 3.5 in \cite{FonsecaMora:2018} shows that indeed $X \in \mathcal{L}(\Psi, C_{\Prob}(T,\R))$. 
Then, for a any given $\epsilon >0$ and by following similar arguments to those used in the proof of Lemma 3.7 in \cite{FonsecaMora:2018}, we can prove that there exists a $\tau_{HS}$-continuous Hilbertian seminorm $p=p(\epsilon)$ on $\Psi$ such that 
\begin{equation} \label{eqEquicoFuncCaracCylProc}
\Exp \left( \sup_{t \in [0,T]} \abs{1-e^{i X_{t}(\psi)}}  \right) \leq \epsilon + p(\psi), \quad \forall \, \psi \in \Psi. 
\end{equation}
Then, if we take $\psi=S \phi$ in \eqref{eqEquicoFuncCaracCylProc} we obtain 
\begin{equation} \label{eqEquicoFuncCaracSazoTopo}
\Exp \left( \sup_{t \in [0,T]} \abs{1-e^{i X_{t}\circ S(\phi)}}  \right) \leq \epsilon + p(S \phi), \quad  \forall \, \phi \in \Phi.  
\end{equation}   
But as $S$, being a Hilbert-Schmidt operator, is continuous from $(\Phi, \tau_{HS})$ into $\Psi$, and $p$ is a separable continuous Hilbertian seminorm on $\Psi$, then 
$q: \Phi \rightarrow \R$ defined by $q(\phi)=p(S \phi)$ for all $\phi \in \Phi$, is therefore $\tau_{HS}$-continuous. We can hence conclude from \eqref{eqEquicoFuncCaracSazoTopo} that the Fourier transforms 
of the family $X\circ S(\phi)= ( X_{t}\circ S ( \phi) : t \in [0,T] )$ are equicontinuous (at zero) in the Sazonov topology on $\Phi$. The result now follows from Theorem \ref{theoRegularizationTheoremLCS}. 
\end{prf}

\begin{coro} \label{coroReguSingHilbSchmUltrab} If $\Psi$ is additionally an  ultrabornological space, the conclusion of Theorem \ref{theoRadonSingleHilSchOper} remain valid if instead of assuming the condition \emph{(2)} we assume that for each $t \geq 0$ the Fourier transform of $X_{t}$ is continuous (at zero) in $\Psi$. 
\end{coro}
\begin{prf} First, if the Fourier transform of $X_{t}$ is continuous (at zero) in $\Psi$ then the mapping $X_{t}: \Psi \rightarrow L^{0} \ProbSpace$ is continuous. Then, since $\Psi$ is ultrabornological, and because for each $\psi \in \Psi$ the real-valued process $X(\psi)$ is continuous (respectively c\`{a}dl\`{a}g), if follows from Proposition 3.10 in \cite{FonsecaMora:2018} that for each $T>0$ the linear mapping  $\psi \mapsto X(\psi)=(X_{t}(\psi): t \in [0,T])$ is continuous from $\Psi$ into $C_{\Prob}(T,\R)$ (respectively into $D_{\Prob}(T,\R)$). This last implies condition \emph{(2)} in Theorem \ref{theoRadonSingleHilSchOper}, and from there the result follows.  
\end{prf}

\begin{exam} \label{examCyliMartingalesHilbSchOperat} Suposse $\Psi$ is equipped with a multi-Hilbertian topology and that $M=(M_{t}: t \geq 0)$ is a cylindrical martingale in $\Psi'$ such that for every $t \geq 0$, the Fourier transform of $M_{t}$ is continuous (at the origin) in $\Psi$. As in the proof of Theorem 5.2 in \cite{FonsecaMora:2018}, we can show that $M$ satisfies conditions $(1)$ and $(2)$ in Theorem \ref{theoRadonSingleHilSchOper}. If $S$ is  a Hilbert-Schmidt operator from $\Phi$ into $\Psi$, then Theorem \ref{theoRadonSingleHilSchOper} shows that $M\circ S= ( M_{t}\circ S : t \geq 0)$ has a $\Phi'$-valued, regular, c\`{a}dl\`{a}g version $Y= (Y_{t}: t \geq 0)$ that is cylindrically a martingale (or a $\Phi'$-valued martingale if $\Phi$ is a separable Hilbert space). Moreover, as in Example \ref{examCyliMartingales}, if $M$ is $r$-integrable for $r \geq 2$, then $Y$ can be constructed such that for any $T>0$, there exists a separable continuous Hilbertian seminorm  $q_{T}$ on $\Phi$ such that $(Y_{t}: t \in [0,T])$ is a $\Phi'_{q_{T}}$-valued c\`{a}dl\`{a}g martingale satisfying 
$\Exp \left( \sup_{t \in [0,T]} q_{T}'(Y_{t})^{r} \right) < \infty.$ 
\end{exam}

\section{Regularization of Cylindrical L\'{e}vy Processes} \label{sectCyliLevy}

We start with our definition of L\'{e}vy processes on the dual of a locally convex space $\Phi$. 

\begin{defi} \label{defiLevyProcess}
A $\Phi'$-valued process $L=( L_{t} :t\geq 0)$ is called a \emph{L\'{e}vy process} if \begin{enumerate}
\item  $L_{0}=0$ a.s., 
\item $L$ has \emph{independent increments}, i.e. for any $n \in \N$, $0 \leq t_{1}< t_{2} < \dots < t_{n} < \infty$ the $\Phi'$-valued random variables $L_{t_{1}},L_{t_{2}}-L_{t_{1}}, \dots, L_{t_{n}}-L_{t_{n-1}}$ are independent,  
\item L has \emph{stationary increments}, i.e. for any $0 \leq s \leq t$, $L_{t}-L_{s}$ and $L_{t-s}$ are identically distributed, and  
\item For every $t \geq 0$ the distribution $\mu_{t}$ of $L_{t}$ is a Radon measure and the mapping $t \mapsto \mu_{t}$ from $[0, \infty)$ into $\goth{M}_{R}^{1}(\Phi')$ is continuous at $0$ in the weak topology.
\end{enumerate}
\end{defi}

Following the definition given in \cite{ApplebaumRiedle:2010} in the context of Banach spaces and in \cite{FonsecaMora:Levy} for duals of nuclear spaces, we introduce the following definition. 

\begin{defi} A cylindrical process $L=( L_{t} :t\geq 0)$ in $\Phi'$ is said to be a \emph{cylindrical L\'{e}vy process} if  $\forall \, n\in \N$, $\phi_{1}, \dots, \phi_{n} \in \Phi$, the $\R^{n}$-valued process $( (L_{t}(\phi_{1}), \dots, L_{t}(\phi_{n})) : t \geq 0)$ is a L\'{e}vy process.  
\end{defi}

One can easily check that the cylindrical process induced by a $\Phi'$-valued L\'{e}vy process is a cylindrical L\'{e}vy process in $\Phi'$ (see \cite{FonsecaMora:Levy}, Lemma 3.7). The converse is in general not true, but sufficient conditions are given in the next result.  
		
\begin{theo} \label{theoCylindrLevyProcessHaveLevyCadlagVersion}
Let $L=( L_{t} :t\geq 0)$  be a cylindrical L\'{e}vy process in $\Phi'$ such that for every $T > 0$, the Fourier transforms of the family $( L_{t}: t \in [0,T] )$ are equicontinuous (at the origin) in the Sazonov topology $\tau_{HS}$ on $\Phi$.  Then, there exists a countably Hilbertian topology $\vartheta$ on $\Phi$ and a $(\widetilde{\Phi_{\vartheta}})'$-valued c\`{a}dl\`{a}g process $Y= ( Y_{t} : t \geq 0)$, such that for every $\phi \in \Phi$, $\inner{Y}{\phi}=( \inner{Y_{t}}{\phi} : t \geq 0)$ is a version of $L(\phi)= ( L_{t}(\phi) : t \geq 0)$. Moreover, $Y$ is a $\Phi'$-valued, regular, c\`{a}dl\`{a}g L\'{e}vy process that is a version of $L$ and that is unique up to indistinguishable versions. 
\end{theo}
\begin{prf}  
The proof is very similar to that of Theorem 3.8 in \cite{FonsecaMora:Levy}. But for completeness we sketch the main arguments. 

First, the existence of the topology $\vartheta$ and the $(\widetilde{\Phi_{\vartheta}})'$-valued c\`{a}dl\`{a}g process $Y= ( Y_{t} : t \geq 0)$ that is a version of $L$ is a consequence of Theorem \ref{theoRegularizationTheoremLCS}. This because each real-valued L\'{e}vy process has a c\`{a}dl\`{a}g version. 

Because for every $\phi_{1}, \dots, \phi_{n} \in \Phi$, $( (\inner{Y_{t}}{\phi_{1}}, \dots, \inner{Y_{t}}{\phi_{n}}: t \geq 0)$ is a $\R^{n}$-valued L\'{e}vy process, and since $Y$ is a regular $\Phi'$-valued process, then $Y=0$ $\Prob$-a.e. and $Y$ has independent and stationary increments (see Propositions 2.2 and 2.3 in \cite{FonsecaMora:Levy}).

Now, since $(\widetilde{\Phi_{\vartheta}})'$ is a Souslin space and $Y$ is $(\widetilde{\Phi_{\vartheta}})'$-valued, for each $t \geq 0$ the probability distribution $\mu_{t}$ of $Y_{t}$ is a Radon measure on $(\widetilde{\Phi_{\vartheta}})'$ (see \cite{BogachevMT}, Theorem 7.4.3, p.85). Moreover, as the canonical inclusion from $(\widetilde{\Phi_{\vartheta}})'$ into $\Phi'$ is linear and continuous, then $\mu_{t}$ is also a Radon measure on $\Phi'$. 

Finally, given $0 \leq t \leq T$, let  $\{ s_{\alpha} \} \subseteq [0,T]$ be a net satisfying $\lim_{\alpha} s_{\alpha}=t$. For any given $\phi \in \Phi$, $(\inner{Y_{t}}{\phi}: t \geq 0) $ is a L\'{e}vy process and hence $t \mapsto \mu_{t} \circ \pi_{\phi}^{-1}$ is continuous in the weak topology. Then,  $\lim_{\alpha} \widehat{\mu}_{s_{\alpha}}(\phi)=\widehat{\mu}_{t}(\phi)$ for all $\phi \in \Phi$. Moreover, our hypothesis on $L$ and the fact that $Y$ is a version of $L$ imply that the family $( \widehat{\mu}_{t} : t \in [0,T])$ is equicontinuous in the Sazonov topology $\tau_{HS}$, and hence the family $( \mu_{r}: r \in [0,T])$ is uniformly tight (see \cite{DaleckyFomin}, Lemma III.2.3, p.103-4). Then, Prokhorov's theorem shows that $( \mu_{s_{\alpha}} )$ is relatively compact, and then we conclude that $\lim_{\alpha} \mu_{s_{\alpha}}=\mu_{t}$ in the weak topology (see \cite{VakhaniaTarieladzeChobanyan}, Theorem IV.3.1, p.224-5). We conclude that the map $t \mapsto \mu_{t}$ is continuous in the weak topology. 
\end{prf}

The following result shows that we can also regularize a cylindrical L\'{e}vy process through a Hilbert-Schmidt operator. 

\begin{theo}  \label{theoRadonLevyBySingleHilSchOper}
Let $\Psi$ be a locally convex space equipped with a multi-Hilbertian topology and let $L=(L_{t}: t\geq 0)$  be a cylindrical L\'{e}vy process in $\Psi'$ such that for every $T > 0$, the Fourier transforms of the family $( L_{t}: t \in [0,T] )$ are equicontinuous (at the origin) on $\Psi$. 
Let $S$ be a Hilbert-Schmidt operator from $\Phi$ into $\Psi$. 
Then, there exists a countably Hilbertian topology $\vartheta$ on $\Phi$ and a $(\widetilde{\Phi_{\vartheta}})'$-valued c\`{a}dl\`{a}g process $Y= ( Y_{t} : t \geq 0)$, such that for every $\phi \in \Phi$, $\inner{Y}{\phi}=( \inner{Y_{t}}{\phi} : t \geq 0)$ is a version of $L \circ S(\phi)= ( L_{t}\circ S ( \phi) : t \geq 0)$. Moreover, $Y$ is a $\Phi'$-valued, regular, c\`{a}dl\`{a}g L\'{e}vy process that is a version of $L \circ S$ (unique up to indistinguishable versions). 
\end{theo}
\begin{prf}
It is clear that $L \circ S(\phi)= ( L_{t}\circ S ( \phi) : t \geq 0)$ is a cylindrical L\'{e}vy process in $\Phi'$. Moreover, from the arguments given in the proof of Theorem \ref{theoRadonSingleHilSchOper}, we conclude that for each $T>0$ the Fourier transforms of the family $L \circ S(\phi)= ( L_{t}\circ S ( \phi) : t \in [0,T])$ are equicontinuous (at the origin) in the Sazonov topology $\tau_{HS}$ on $\Phi$. The result now follows from Theorem \ref{theoRadonLevyBySingleHilSchOper}. 
\end{prf}

By following similar arguments to those used in the proof of Corollary \ref{coroReguSingHilbSchmUltrab} together with Theorem \ref{theoRadonLevyBySingleHilSchOper} we conclude the following:
  
\begin{coro} \label{coroReguLevyBySingHilbSchmUltrab} If $\Psi$ is additionally an  ultrabornological space, the conclusion of Theorem \ref{theoRadonLevyBySingleHilSchOper} remain valid if  we only assume that for each $t \geq 0$ the Fourier transform of $X_{t}$ is continuous (at zero) in $\Psi$. 
\end{coro}

\begin{exam} \label{examRadoSingHSOper}
Let $\Phi$ and $\Psi$ denote two Hilbert spaces, $\Psi$ being separable. Let $S$ be a  Hilbert-Schmidt operator from $\Phi$ into $\Psi$, and $X=(X_{t}: t \geq 0)$ be a cylindrical process in $\Psi$ such that for each  $\psi \in \Psi$,  $X(\psi)=( X_{t}(\psi) : t \geq 0)$ has a c\`{a}dl\`{a}g version and such that for each $t \geq 0$ the the Fourier transform of $X_{t}$ is continuous (at the origin) on $\Psi$. It is a direct consequence of Corollary \ref{coroReguSingHilbSchmUltrab} that $X \circ S=( X_{t} \circ S: t \geq 0)$ has a $\Phi$-valued c\`{a}dl\`{a}g version $Y=(Y_{t}: t \geq 0)$. If $X$ is a cylindrical L\'{e}vy then Corollary \ref{coroReguLevyBySingHilbSchmUltrab} shows that indeed $Y$ is a L\'{e}vy process. 
\end{exam}

\begin{exam} \label{examReguGelfandTriplet} Suppose we have a nuclear space $\Phi$, a separable Hilbert space $\Psi$, and a continuous linear operator $S: \Phi \rightarrow \Psi$ (e.g. $S$ could be an embedding). Suppose that $X=(X_{t}: t \geq 0)$ is a cylindrical process in $\Psi$ such that for each  $\psi \in \Psi$,  $X(\psi)=( X_{t}(\psi) : t \geq 0)$ has a c\`{a}dl\`{a}g version and such that for each $t \geq 0$ the the Fourier transform of $X_{t}$ is continuous (at the origin) on $\Psi$. 

Since the nuclear topology on $\Phi$ coincides with its Sazonov topology, then $S$ is  a Hilbert-Schmidt operator (see Section \ref{sectSazonovTopo}). Then, Corollary \ref{coroReguSingHilbSchmUltrab} shows that $X \circ S=( X_{t} \circ S: t \geq 0)$ has a $\Phi'$-valued c\`{a}dl\`{a}g version $Y=(Y_{t}: t \geq 0)$. If $S$ is an embedding, and $\norm{\cdot}$ denotes the Hilbertian norm on $H$, then $p(\phi)=\norm{S(\phi)}$, $\phi \in \Phi$, is a continuous Hilbertian seminorm on $\Phi$ such that $X \circ S$ is $p$-continuous. So indeed, Corollary \ref{coroRegulSingleSeminorm} shows that there exists a continuous Hilbertian seminorm $q$ on $\Phi$, $p \leq q$, such that $i_{p,q}$ is Hilbert-Schmidt and $Y$ is a $\Phi'_{q}$-valued  c\`{a}dl\`{a}g version of $X \circ S$. 
Moreover, if $X$ is a cylindrical L\'{e}vy process then Corollary \ref{coroReguLevyBySingHilbSchmUltrab} shows that indeed $Y$ is a L\'{e}vy process. 

A practical case that fits into this situation is when we take for example $\Phi= \mathcal{S}(\R^{d})$ and $\Psi=L^{2}(\R^{d})$, and $X$ a cylindrical process in $L^{2}(\R^{d})$ satisfying the conditions given above. Then, because the well-known fact that the canonical embedding $I:\mathcal{S}(\R^{d}) \rightarrow L^{2}(\R^{d})$ is linear and  continuous, then $X \circ I$ has a $\mathcal{S}'(\R^{d})$-valued c\`{a}dl\`{a}g version (which is a L\'{e}vy process if $X$ is cylindrical L\'{e}vy). 
\end{exam}



\section{Discussion and Comparison With Known Results} \label{sectDiscuComp}


To the extend of our knowledge Theorem \ref{theoRegularizationTheoremLCS} is the first attempt to provide sufficient conditions for the existence of continuous and c\`{a}dl\`{a}g versions to cylindrical processes in general locally convex spaces. However, for particular classes of locally convex spaces there are some other interesting works for  which Theorem \ref{theoRegularizationTheoremLCS} constitutes a generalization, in particular we have:
\begin{enumerate}
\item Theorem \ref{theoRegularizationTheoremLCS} extends the conclusions of It\^{o} and Nawata regularization theorem from multi-Hilbertian spaces to the context of general locally convex spaces (see \cite{ItoNawata:1983} and Theorem 2.3.2 in \cite{Ito}).

\item Since every cylindrical probability measure on $\Phi'$ possesses a canonical cylindrical random variable in $\Phi'$ defined on some probability space (see Section II.V.2 in \cite{SchwartzRM}, p.256-8), then Theorem \ref{theoRegularizationTheoremLCS} offers an alternative proof to the Minlos-Bochner theorem (Theorem III.1.1 in \cite{DaleckyFomin}).  

\item  Theorem \ref{theoRegularizationTheoremLCS} includes the regularization theorems for the existence of continuous and c\`{a}dl\`{a}g versions to $\Phi'$-valued processes when $\Phi$ is \begin{inparaenum}[(i)] \item  a nuclear Fr\'{e}chet space (Mitoma \cite{Mitoma:1983}), \item a countable inductive limit of nuclear Fr\'{e}chet spaces (Fouque \cite{Fouque:1984} and Fernique \cite{Fernique:1989}), \item or a separable nuclear space (see Martias \cite{Martias:1988}). \end{inparaenum}
To obtain these conclusions we use the fact that $\tau_{HS}=\tau$ when $(\Phi,\tau)$ is nuclear, together with Theorem 2.10 and Proposition 3.10 in \cite{FonsecaMora:2018}. 

\item When $\Phi$ is a (general) nuclear space, Theorem \ref{theoRegularizationTheoremLCS} directly implies the regularization theorem for cylindrical processes in $\Phi'$  (\cite{FonsecaMora:2018}, Theorem 3.2). Again, in this case to get this conclusion we use the fact that $\tau_{HS}=\tau$. In a similar way, our Corollary 3.5 implies Theorem 4.1 in \cite{FonsecaMora:2018}.   

\end{enumerate}

The version of the regularization theorem for cylindrical processes which possesses finite moments given in Therem \ref{theoRegulFiniteMomentsLCS} and the conclusions obtained in Example \ref{examCyliMartingales} for cylindrical martingales, generalize the results obtained in Theorems 4.2 and 5.2 in \cite{FonsecaMora:2018} (there proved under the assumption thet $\Phi$ is a nuclear space), and also the results for martingales in duals of nuclear Fr\'{e}chet spaces in \cite{Mitoma:1981}.

In a similar way, Theorem \ref{theoRadonSingleHilSchOper} includes some other results from the literature. In particular, the result in Example \ref{examRadoSingHSOper}, often known as \emph{radonification by a single Hilbert-Schmidt operator}, was proven first by Badrikian and \"{U}st\"{u}nel in \cite{BadrikianUstunel:1996} and later by Jakubowski et al. in \cite{JakubowskiEtAl:2002}. In both works their motivation was to show the radonification of a cylindrical semimartingale through three Hilbert-Schmidt operators in \cite{BadrikianUstunel:1996} and by one Hilbert-Schmidt operator in \cite{JakubowskiEtAl:2002}. 
The radonification of a cylindrical L\'{e}vy process through a single Hilbert-Schmidt operator is somewhat a known fact but we have no knowledge of any reference with a formal proof of it. Its extension to general locally convex spaces as in Theorem \ref{theoRadonLevyBySingleHilSchOper} as well the conclusion in Example \ref{examReguGelfandTriplet} are new. 

If $(\Phi, \tau)$ is a nuclear space,  Theorem \ref{theoCylindrLevyProcessHaveLevyCadlagVersion} coincides with Theorem 3.8 in \cite{FonsecaMora:Levy}. However, outside the nuclear space setting, to the extend of our knowledge the conclusions of Theorem \ref{theoCylindrLevyProcessHaveLevyCadlagVersion} are new, even in the Hilbert space setting.  

Appart from the Sazonov topology, other topologies has been considered on a locally convex space to provide sufficient conditions for a cylindrical measure to extend to a Radon measure. One of these topologies is the following (see e.g. \cite{VakhaniaTarieladzeChobanyan}). Let $\Psi$ denote a Hilbert space. A symmetric, positive, nuclear operator in $\Psi$ is called an $S$-operator. 

Consider a locally convex space $\Phi$. Denote by $\overline{S}(\Phi, \Phi')$ the class of operators $R: \Phi \rightarrow \Phi'$ of the form $R= v' \circ S \circ v$, where $S$ is an $S$-operator in a separable Hilbert space $\Psi$ and $v: \Phi \rightarrow \Psi$ is a continuous linear operator. Denote by $\tau_{\overline{S}}=\tau_{\overline{S}}(\Phi,\Phi')$ the weakest vector topology in $\Phi$ with respect to which all the quadratic forms $\phi \mapsto \inner{R \phi}{\phi}$, $R \in \overline{S}(\Phi,\Phi')$ are continuous. 
If $\Phi$ is a separable Hilbert space, then it is easy to check that $\tau_{HS}=\tau_{\overline{S}}$, but in general, e.g. if $\Phi$ is Banach, $\tau_{\overline{S}}$ is finer than $\tau_{HS}$. From this last fact in combination with Theorem \ref{theoRegularizationTheoremLCS} we obtain the following: 

\begin{coro}\label{coroRegulTheoLCSForTauSTopo}
Suppose that $X$ is a cylindrical process in $\Phi'$ such that for each $\phi \in \Phi$, the real-valued process $X(\phi)$ has a continuous (respectively c\`{a}dl\`{a}g) version. The conclusions of Theorem \ref{theoRegularizationTheoremLCS} remain valid if we assume that for each $T>0$, the the Fourier transforms 
of the family $( X_{t}: t \in [0,T] )$ are equicontinuous (at the origin) in the topology $\tau_{\overline{S}}$. 
\end{coro}

If again we consider a cylindrical probability measure on $\Phi'$, and consider its  canonical cylindrical random variable in $\Phi'$, then Corollary \ref{coroRegulTheoLCSForTauSTopo} can be used to show the version of the Minlos-Sazonov theorem based on continuity on the $\tau_{\overline{S}}$-topology given in Theorem IV.4.1. in \cite{VakhaniaTarieladzeChobanyan}. 

\textbf{Acknowledgements} { 
The author is grateful to Dar\'{i}o Mena-Arias for stimulating discussions. The author acknowledge The University of Costa Rica for providing financial support through the grant ``Pry01-1692-2019-An\'{a}lisis Estoc\'{a}stico con Procesos Cil\'{i}ndricos''. 
}

\end{document}